\documentclass[11pt]{article}

\usepackage{amssymb}
\usepackage{amsmath}
\usepackage{amscd}

\newtheorem{theorem}{Theorem}[section]

\begin{document}

\title{On Sojourn Times in the Finite Capacity $M/M/1$ Queue with Processor Sharing}
\author{
Qiang Zhen\thanks{
Department of Mathematics, Statistics, and Computer Science,
University of Illinois at Chicago, 851 South Morgan (M/C 249),
Chicago, IL 60607-7045, USA.
{\em Email:} qzhen2@uic.edu.}
\and
and
\and
Charles Knessl\thanks{
Department of Mathematics, Statistics, and Computer Science,
University of Illinois at Chicago, 851 South Morgan (M/C 249),
Chicago, IL 60607-7045, USA.
{\em Email:} knessl@uic.edu.\
\newline\indent\indent{\bf Acknowledgement:} This work was partly supported by NSF grant DMS 05-03745 and NSA grant H 98230-08-1-0102.
}}
\date{ }
\maketitle

\begin{abstract}
We consider a processor shared $M/M/1$ queue that can accommodate at most a finite number $K$ of customers. We give an exact expression for the sojourn time distribution in the finite capacity model, in terms of a Laplace transform. We then give the tail behavior, for the limit $K\to\infty$, by locating the dominant singularity of the Laplace transform.

{\bf keywords:} finite capacity, processor sharing, sojourn time, tail behavior
\end{abstract}

\section{Introduction}
Processor sharing (PS) is one of the most interesting service disciplines in queueing theory. The $M/M/1$-PS queue assumes Poisson arrivals with rate $\lambda$ and exponential i.i.d. service times with mean $1/\mu$. The traffic intensity is $\rho=\lambda/\mu$. Here we consider a system which can accommodate at most $K$ customers. If the system is filled to capacity, we assume that further arrivals are turned away and lost.

Most past work on PS models deals with systems that have an infinite capacity of customers. In \cite{COF}, Coffman, Muntz, and Trotter analyzed the $M/M/1$-PS queue, and derived an expression for the Laplace transform of the sojourn time distribution, conditioned on both the number of other customers seen, and the amount of service required, by an arriving customer. Using the results in \cite{COF}, Morrison \cite{MO} obtains asymptotic expansions for the unconditional sojourn time distribution in the heavy traffic limit, where the Poisson arrival rate $\lambda$ is nearly equal to the service rate $\mu$ (thus $\rho=\lambda/\mu\uparrow 1$).

A service discipline seemingly unrelated to PS is that of random order service (ROS), where customers are chosen for service at random. In \cite{PO} Pollaczek derives an explicit integral representation for the generating function of the waiting time distribution $\mathbf{W}_{\mathrm{ROS}}$, from which the following tail behavior is obtained
\begin{equation}\label{S1_tail}
\Pr\left[\mathbf{W}_{\mathrm{ROS}}>t\right]\sim e^{-\alpha t-\beta t^{1/3}}\gamma t^{-5/6}, t\rightarrow\infty.
\end{equation}
Here $\alpha$, $\beta$ and $\gamma$ are explicitly computed constants, with $\alpha=(1-\sqrt{\rho})^2$. Cohen \cite{COH} establishes the following relationship between the sojourn time $\mathbf{V}_{\mathrm{PS}}$ in the PS model and the waiting time $\mathbf{W}_{\mathrm{ROS}}$ in the ROS model,
\begin{equation}\label{S1_equi}
\rho\Pr\left[\mathbf{V}_{\mathrm{PS}}>t\right]=\Pr[\mathbf{W}_{\mathrm{ROS}}>t],
\end{equation}
which extends also to the more general $G/M/1$ case. In \cite{BO} relations of the form (\ref{S1_equi}) are explored for other models, such as finite capacity queues, repairman problems, and networks.

The present finite capacity model will have purely exponential behavior for $t\to\infty$, and thus the subexponential ($e^{-\beta\,t^{1/3}}$) and algebraic ($t^{-5/6}$) factors that appear in (\ref{S1_tail}) will be absent. Writing $\Pr[\mathbf{V}>t]\sim C\,e^{-\delta t}$ for $t\to\infty$ for the finite capacity model, the relaxation rate $\delta=\delta(K)$ depends on the capacity, and we shall study its behavior for $K\to\infty$ and for various values of $\rho$. In particular, if $\rho<1$ and $K\to\infty$ it will prove instructive to compare the tail behaviors when $K=\infty$ (where (\ref{S1_tail}), (\ref{S1_equi}) apply) and when $K$ is large but finite.

Knessl \cite{knessl1990} obtains asymptotic expansions for the first two conditional sojourn time moments for the finite capacity model. The assumption is that $K\gg 1$ and separate analyses are carried out for the cases $\rho<1$, $\rho=1$ and $\rho>1$. In \cite{knessl1993} he obtains expansions for the sojourn time distribution by performing asymptotic analyses for the three scales $\rho-1=O(K^{-1})$; $\rho-1=bK^{-1/2},\, b>0$ and $\rho>1$.

The main purpose of this note is to give an explicit exact expression for the Laplace transform of the sojourn time distribution in the finite capacity $M/M/1$-PS model. We apply a Laplace transform to the basic evolution equation and solve it using a discrete Green's function. The solution is summarized in Theorem 2.1. Then in Theorem 2.2 we locate the dominant singularity, which leads to the tail behavior, for $K\gg 1$.

\section{Summary of results}
We set the service rate $\mu=1$. Then the traffic intensity is $\rho=\lambda>0$. We define the conditional density of the sojourn time $\mathbf{V}$ by
\begin{equation*}\label{S2_def_pnt}
p_n(t)dt=\Pr\big[\mathbf{V}\in(t,t+dt)\big|\mathbf{N(0^-)}=n], \quad 0\leq n\leq K-1.
\end{equation*}
Here $\mathbf{N(0^-)}$ is the number of customers present in the system immediately before the tagged customer arrives.  From \cite{knessl1993}, the quantity $p_n(t)$ satisfies the evolution equation
\begin{equation}\label{S2_recu}
p'_n(t)=\rho\; p_{n+1}(t)-(1+\rho)\; p_n(t)+\frac{n}{n+1}\;p_{n-1}(t),\quad 0\leq n\leq K-1,
\end{equation}
for $t>0$, with the initial condition $p_n(0)=\frac{1}{n+1}$. If $n=K-1$ the term $\rho\, p_{n+1}(t)$ is absent, and we shall use the ``artificial" boundary condition
$p_{_{K-1}}(t)=p_{_K}(t)$.
Taking the Laplace transform of (\ref{S2_recu}) with $\widehat{p}_n(\theta)=\int^\infty_0 p_n(t)e^{-\theta t}dt$ and multiplying by $n+1$, we have
\begin{equation}\label{S2_recu2}
(n+1)\;\rho\;\widehat{p}_{n+1}(\theta)-(n+1)\;(1+\rho+\theta)\;\widehat{p}_n(\theta)+n\;\widehat{p}_{n-1}(\theta)=-1
\end{equation}
for $n=0,1,...,K-1$ with the boundary condition
\begin{equation}\label{S2_bc2}
\widehat{p}_{_{K-1}}(\theta)=\widehat{p}_{_K}(\theta).
\end{equation}

Solving the recurrence equation (\ref{S2_recu2}) with (\ref{S2_bc2}), we obtain the following result.

\begin{theorem} \label{th1}
The Laplace-Stieltjes transform of the conditional sojourn time density has the following form:
\begin{equation}\label{S2_th1_phat}
\widehat{p}_n(\theta)=MG_n\sum_{l=0}^n\rho^lH_l+MH_n\sum_{l=n+1}^{K-1}\rho^lG_l-\frac{\Delta G_K}{\Delta H_K}MH_n\sum_{l=0}^{K-1}\rho^lH_l
\end{equation}
for $n=0,1,...,K-2$ and
\begin{equation}\label{S2_th1_phat_K}
\widehat{p}_{_{K-1}}(\theta)=\frac{1}{K\rho^K\Delta H_K}\sum_{l=0}^{K-1}\rho^lH_l.
\end{equation}
Here
\begin{equation}\label{S2_th1_M}
M=M(\theta)\equiv z_-\Big(\frac{z_+}{z_-}\Big)^\alpha,
\end{equation}
\begin{equation}\label{S2_th1_G}
G_n=G_n(\theta)\equiv\int_0^{z_-}z^n(z_+-z)^{-\alpha}(z_--z)^{\alpha-1}dz,
\end{equation}
\begin{equation}\label{S2_th1_H}
H_n=H_n(\theta)\equiv\frac{e^{i\alpha\pi}}{2\pi i}\int_{\mathcal{C}}z^n(z_+-z)^{-\alpha}(z-z_-)^{\alpha-1}dz,
\end{equation}
$\mathcal{C}$ is a closed contour in the complex $z$-plane that encircles the segment $[z_-,z_+]$ of the real axis counterclockwise and
\begin{equation*}\label{S2_th1_deltag}
\Delta G_K=G_K-G_{K-1},
\end{equation*}
\begin{equation}\label{S2_th1_deltah}
\Delta H_K=H_K-H_{K-1},
\end{equation}
\begin{equation}\label{S2_th1_z}
z_\pm=z_\pm(\theta)\equiv\frac{1}{2\rho}\Big[1+\rho+\theta\pm\sqrt{(1+\rho+\theta)^2-4\rho}\Big],
\end{equation}
\begin{equation}\label{S2_th1_alpha}
\alpha=\alpha(\theta)\equiv\frac{z_+}{z_+-z_-}.
\end{equation}

\end{theorem}

The singularities of $\widehat{p}_n(\theta)$ are poles, which are all real, and solutions of $H_K(\theta)=H_{K-1}(\theta)$. A spectral expansion can be given in terms of these poles, but we believe that (\ref{S2_th1_phat}) is more useful for certain asymptotic analyses, such as $K\to\infty$. Below we give the least negative pole $\theta_s$, which is also the tail exponent, i.e. $\lim_{t\to\infty}\{t^{-1}\log\Pr[\mathbf{V}>t]\}$.

\begin{theorem} \label{th2}
The dominant singularity $\theta_s$ in the $\theta$-plane has the following asymptotic expansions, for $K\to\infty$:

\begin{enumerate}

\item $\rho<1$,
\begin{equation}\label{S2_th2_rho<1}
\theta_s=-(1-\sqrt{\rho})^2-\frac{\sqrt{\rho}}{K}+\frac{\sqrt{\rho}\,r_0}{K^{4/3}}-\frac{8\sqrt{\rho}\,r_0^2}{15K^{5/3}}+O(K^{-2}),
\end{equation}
where $r_0\approx -2.3381$ is the largest root of the Airy function $Ai(z)$.

\item $\rho=1+\eta\,K^{-2/3}$ with $\eta=O(1)$,
\begin{equation}\label{S2_th2_rho=1}
\theta_s=-\frac{1}{K}+\frac{r_1}{K^{4/3}}-\frac{16r_1^3+8\eta^2\,r_1^2+(\eta^4+19\eta)r_1+\eta^3+9}{30r_1\,K^{5/3}}+O(K^{-2}),
\end{equation}
where $r_1=r_1(\eta)$ is the largest root of
\begin{equation}\label{S2_th2_rho=1_r1}
\frac{Ai'(r_1+\eta^2/4)}{Ai(r_1+\eta^2/4)}=-\frac{\eta}{2}.
\end{equation}

\item $\rho>1$,
\begin{equation}\label{S2_th2_rho>1}
\theta_s=-\frac{1}{K}-\frac{1}{(\rho-1)K^{2}}-\frac{1}{(\rho-1)^2K^{3}}+\frac{\rho^2+1}{(\rho-1)^4K^4}+O(K^{-5}).
\end{equation}

\end{enumerate}

\end{theorem}

Comparing $e^{\theta_st}$ with (\ref{S2_th2_rho<1}) to (\ref{S1_tail}) we see that they both contain the dominant factor $e^{-(1-\sqrt{\rho})^2t}$ but whereas (\ref{S1_tail}) has the subexponential factor $e^{-\beta\,t^{1/3}}$, (\ref{S2_th2_rho<1}) leads to purely exponential correction terms that involve the maximal root of the Airy function.

\section{Brief derivations}

We use a discrete Green's function to derive (\ref{S2_th1_phat}) and (\ref{S2_th1_phat_K}). Consider the recurrence equation (\ref{S2_recu2}) and (\ref{S2_bc2}). The discrete Green's function $\mathcal{G}(\theta;n,l)$ satisfies
\begin{eqnarray}\label{S3_dgf}
&&(n+1)\rho\,\mathcal{G}(\theta;n+1,l)-(n+1)(1+\rho+\theta)\,\mathcal{G}(\theta;n,l)\nonumber\\
&&\quad\quad\quad+n\,\mathcal{G}(\theta;n-1,l)=-\delta(n,l),\quad (n,l=0,1,...,K-1)
\end{eqnarray}
and
\begin{equation}\label{S3_dgf_bc}
\mathcal{G}(\theta;K,l)=\mathcal{G}(\theta;K-1,l),\quad (l=0,1,...,K),
\end{equation}
where $\delta(n,l)=1_{\{n=l\}}$ is the Kronecker delta.
To construct the Green's function requires that we have two linearly independent solutions to
\begin{equation}\label{S3_dgf_Homo}
(n+1)\rho\,G(\theta;n+1,l)-(n+1)(1+\rho+\theta)\,G(\theta;n,l)+n\,G(\theta;n-1,l)=0,
\end{equation}
which is the homogeneous version of (\ref{S3_dgf}). We seek solutions of (\ref{S3_dgf_Homo}) in the form
$ G_n=\int_\mathcal{D}z^ng(z)dz,$
where the function $g(z)$ and the path $\mathcal{D}$ of integration in the complex $z$-plane are to be determined. Using the above form in (\ref{S3_dgf_Homo}) and integrating by parts yields
\begin{eqnarray}\label{S3_IBP}
&&z^ng(z)\big[\rho z^2-(1+\rho+\theta)z+1\big]\Big|_\mathcal{D}\nonumber\\
&&\quad\quad -\int_\mathcal{D}z^n\big[(\rho z^2-(1+\rho+\theta)z+1)g'(z)+\rho zg(z)\big]dz=0.
\end{eqnarray}
The first term represents contributions from the endpoints of the contour $\mathcal{D}$.

If (\ref{S3_IBP}) is to hold for all $n=0,1,...,K-1$ the integrand must vanish, so that $g(z)$ must satisfy the differential equation
$$\big[\rho z^2-(1+\rho+\theta)z+1\big]g'(z)+\rho zg(z)=0.$$
We denote the roots of $\rho z^2-(1+\rho+\theta)z+1=0$ by $z_+$ and $z_-$, with $z_+>z_->0$ for real $\theta$. These are given by (\ref{S2_th1_z}) and if $\alpha$ is defined by (\ref{S2_th1_alpha}), the solution for $g(z)$ is
$ g(z)=(z_+-z)^{-\alpha}(z_--z)^{\alpha-1}.$
By appropriate choices of the contour $\mathcal{D}$, we obtain the two linearly independent solutions $G_n$ (cf. (\ref{S2_th1_G})) and $H_n$ (cf. (\ref{S2_th1_H})). We note that $G_n$ decays as $n\to\infty$, and is asymptotically given by
\begin{equation}\label{S3_Gasymp}
G_n\sim\frac{\Gamma(\alpha)}{n^\alpha}z_-^{\alpha+n}(z_+-z_-)^{-\alpha},\quad n\to\infty.
\end{equation}
However, $G_n$ becomes infinite as $n\to-1$, and $nG_{n-1}$ goes to a nonzero limit as $n\to 0$. Thus $G_n$ is not an acceptable solution at $n=0$. $H_n$ is finite as $n\to-1$, but grows as $n\to\infty$, with 
\begin{equation}\label{S3_Hasymp}
H_n\sim\frac{n^{\alpha-1}}{\Gamma(\alpha)}z_+^{n+1-\alpha}(z_+-z_-)^{\alpha-1},\quad n\to\infty.
\end{equation}

Thus, the discrete Green's function can be represented by
\begin{eqnarray}\label{S3_G_step}
\mathcal{G}(\theta;n,l) = \left\{ \begin{array}{ll}
AH_n+BG_n & \textrm{if $l\leq n\leq K$}\\
CH_n & \textrm{if $0\leq n\leq l\leq K$}.
\end{array} \right.
\end{eqnarray}
Here $A$, $B$ and $C$ depend upon $\theta$, $l$ and $K$. Using continuity of $\mathcal{G}$ at $n=l$ and the boundary condition (\ref{S3_dgf_bc}), we obtain from (\ref{S3_G_step})
$$ A=\frac{H_l\Delta G_{{K}}}{H_l\Delta G_{{K}}-G_l\Delta H_{{K}}}C,\quad \quad B=-\frac{H_l\Delta H_{K}}{H_l\Delta G_{{K}}-G_l\Delta H_{{K}}}C.$$
Hence, (\ref{S3_G_step}) can be rewritten as
\begin{eqnarray*}\label{S3_G_step2}
\mathcal{G}(\theta;n,l) = \left\{ \begin{array}{ll}
{\displaystyle\frac{H_n\Delta G_K-G_n\Delta H_K}{H_l\Delta G_{{K}}-G_l\Delta H_{{K}}}\,CH_l} & \textrm{if $l\leq n\leq K$}\bigskip\\
CH_n & \textrm{if $0\leq n\leq l\leq K$}.
\end{array} \right.
\end{eqnarray*}

To determine $C$, we let $n=l$ in (\ref{S3_dgf}) and use the fact that both $G_l$ and $H_l$ satisfy (\ref{S3_dgf_Homo}) with $n=l$, which gives 
\begin{equation}\label{S3_C}
C=\frac{G_l\Delta H_{{K}}-H_l\Delta G_{{K}}}{(l+1)\,\rho\,\Delta H_K(G_l\,H_{l+1}-G_{l+1}\,H_l)}.
\end{equation}
From (\ref{S3_dgf_Homo}) we can infer a simple difference equation for the discrete Wronskian $G_l\,H_{l+1}-G_{l+1}\,H_l$, whose solution we write as
\begin{equation}\label{S3_G1}
G_l\,H_{l+1}-G_{l+1}\,H_l=\frac{1}{C_1\,(l+1)\,\rho\,^l},
\end{equation}
where $C_1=C_1(\theta)$ depends upon $\theta$ only. Letting $l\to\infty$ in (\ref{S3_G1}) and using the asymptotic results in (\ref{S3_Gasymp}) and (\ref{S3_Hasymp}), we determine $C_1$ as $C_1=\rho M$ and then by (\ref{S3_C}) obtain
$$C=\frac{G_l\Delta H_{{K}}-H_l\Delta G_{{K}}}{\Delta H_K}\rho^lM,$$
where $M=M(\theta)$ is defined by (\ref{S2_th1_M}).

Then, we multiply (\ref{S3_dgf}) by the solution $\widehat{p}_l(\theta)$ to (\ref{S2_recu2}) and sum over  $0\leq l\leq K-1$. After some manipulation and using (\ref{S3_dgf_bc}), this yields
$$\widehat{p}_n(\theta)=\sum_{l=0}^{K-1}\mathcal{G}(\theta;n,l),\quad n=0,1,...,K-1,$$
which is equivalent to (\ref{S2_th1_phat}) and (\ref{S2_th1_phat_K}).

In the remainder of this section, we obtain the dominant singularity of $\widehat{p}_n(\theta)$ in the $\theta$-plane. From (\ref{S2_th1_phat}), the dominant singularity comes from solving $\Delta H_{{K}}=H_{{K}}-H_{{K-1}}=0$.

We first consider $\rho<1$. We evaluate $H_n$ in (\ref{S2_th1_H}) by branch cut integration, which yields
\begin{equation}\label{S3_Hn}
H_n=\frac{\sin({\alpha\pi})}{\pi}\int_{z_-}^{z_+}\xi^n\,(\xi-z_-)^{\alpha-1}(z_+-\xi)^{-\alpha}d\xi.
\end{equation}
Thus, we have
\begin{equation}\label{S3_deltah}
\Delta H_K=\frac{\sin({\alpha\pi})}{\pi}\int_{z_-}^{z_+}\frac{\xi-1}{\xi-z_-}\xi^{K-1}\left(\frac{\xi-z_-}{z_+-\xi}\right)^\alpha d\xi.
\end{equation}
Since the factor $\sin(\alpha \pi)$ appears in both the numerators and denominators in (\ref{S2_th1_phat}) and (\ref{S2_th1_phat_K}), zeros of $\sin(\alpha \pi)$ do not correspond to poles of $\widehat{p}_n(\theta)$. To be at a pole we must have the integral in (\ref{S3_deltah}) vanish.
When $\rho<1$ and $K\gg 1$, the effects of finite capacity should become asymptotically negligible and the finite capacity model may be approximated by the corresponding infinite capacity model. Thus, from (\ref{S1_tail}) the dominant singularity $\theta_s$ should be close to $-\alpha=-(1-\sqrt{\rho})^2$. We scale $\theta=-(1-\sqrt{\rho})^2+s/K$ and let
\begin{equation}\label{S3_xi}
\xi=\frac{z_++z_-}{2}+\frac{z_+-z_-}{2}\,w
\end{equation}
in (\ref{S3_deltah}). After some calculation, the integral in (\ref{S3_deltah}) is approximately given by
\begin{equation}\label{S3_delta_rho<1}
\frac{1-\sqrt{\rho}}{\rho^{K/2}}\int_{-1}^1\Big[g(w,s)+O(K^{-1/2})\Big]\,e^{\sqrt{K}\,f(w,s)}\,dw,
\end{equation}
where
\begin{equation*}\label{S3_f}
f(w,s)=\rho^{-1/4}\,w\,\sqrt{s}+\frac{\rho^{1/4}}{2\sqrt{s}}\log\Big(\frac{1+w}{1-w}\Big),
\end{equation*}
\begin{equation*}\label{S3_g}
g(w,s)=\frac{1}{\sqrt{1-w^2}}\exp\Big[\frac{s(1-w^2)}{2\sqrt{\rho}}\Big].
\end{equation*}
By solving $\frac{\partial}{\partial w}f(w,s)=0$, we find that there are two saddle points, at $w_*=\pm\sqrt{1+\sqrt{\rho}/s}$. In order to be in the oscillatory range of (\ref{S3_delta_rho<1}) the two saddles must coalesce, which occurs at $w=0$ when $s=-\sqrt{\rho}$. Then we introduce the scaling $w=uK^{-1/6}$ and $s=\sqrt{\rho}(-1+rK^{-1/3})$ and expand the functions $f(w,s)$ and $g(w,s)$ as
\begin{equation}\label{S3_f_asym}
f(w,s)=-i\,\Big(ru+\frac{u^3}{3}\Big)K^{-1/2}-i\,\Big(\frac{r^2u}{2}+\frac{ru^3}{6}+\frac{u^5}{5}\Big)K^{-5/6}+O\big(K^{-7/6}\big),
\end{equation}
\begin{equation}\label{S3_g_asym}
g(w,s)=e^{-1/2}\Big[1+\Big(\frac{r}{2}+u^2\Big)K^{-1/3}+O(K^{-2/3})\Big].
\end{equation}
By using the definition of the Airy function
\begin{equation*}
Ai(z)=\frac{1}{2\pi}\int_{-\infty}^{\infty}e^{-i(t^3/3+zt)}dt,
\end{equation*}
along with (\ref{S3_f_asym}) and (\ref{S3_g_asym}), the integral in (\ref{S3_delta_rho<1}) can be evaluated as
\begin{equation}\label{S3_delta_rho<1_2}
2\pi e^{-1/2}K^{-1/6}\bigg\{Ai(r)+K^{-1/3}\Big[\frac{2}{15}rAi(r)+\frac{8}{15}r^2Ai'(r)\Big]+O(K^{-2/3})\bigg\}.
\end{equation}
To find the dominant singularity, (\ref{S3_delta_rho<1_2}) must vanish. It follows that $r$ is asymptotic to the maximal root of $Ai$. To get a refined approximation, we expand $r$ as $r=r_0+\alpha_0K^{-1/3}+O(K^{-2/3})$ and calculate $\alpha_0$ from (\ref{S3_delta_rho<1_2}). Then using $\theta=-(1-\sqrt{\rho})^2-\sqrt{\rho}/K+\sqrt{\rho}\,r\,K^{-4/3}$ we obtain (\ref{S2_th2_rho<1}).

Next, we consider $\rho=1+\eta K^{-2/3}$ with $\eta=O(1)$. We now scale $\theta=-1/K+R\,K^{-4/3}$ and $w=u\,K^{-1/6}$, with $R$ and $u=O(1)$. Then the leading order approximation to the integrand in (\ref{S3_deltah}) is 
\begin{equation}\label{S3_integrand_rho=1}
e^{-K^{1/3}\eta/2-1/2}K^{-5/6}\Big[\big(i\,u-\frac{\eta}{2}\big)+O(K^{-1/3})\Big]\exp\Big\{-i\Big[\frac{u^3}{3}+\big(R+\frac{\eta^2}{4}\big)\,u\Big]\Big\}.
\end{equation}
Integrating (\ref{S3_integrand_rho=1}) over $u$ now leads to
\begin{equation}\label{S3_delta_rho=1_2}
2\pi e^{-K^{1/3}\eta/2-1/2}K^{-5/6}\Big[-\frac{\eta}{2}Ai\big(R+\frac{\eta^2}{4}\big)-Ai'\big(R+\frac{\eta^2}{4}\big)+O(K^{-1/3})\Big].
\end{equation}
Setting (\ref{S3_delta_rho=1_2}) equal to zero we obtain the expansion $R=r_1+\alpha_1K^{-1/3}+O(K^{-2/3})$, which leads to $\theta_s\sim-1/K+r_1K^{-4/3}$, where $r_1$ satisfies (\ref{S2_th2_rho=1_r1}). We note that to obtain the $O(K^{-5/3})$ term in (\ref{S2_th2_rho=1}), we needed the $O(K^{-1/3})$ correction terms in (\ref{S3_integrand_rho=1}) and (\ref{S3_delta_rho=1_2}), which we do not give here.

Finally, we consider $\rho>1$. We use (\ref{S2_th1_H}) and rewrite $\Delta H_K$ as
\begin{equation}\label{S3_delta_rho>1}
\Delta H_K=\frac{e^{i\alpha\pi}}{2\pi i}\int_{\mathcal{C}}z^{K-1}(z-1)\frac{(z-z_-)^{\alpha-1}}{(z_+-z)^\alpha}dz.
\end{equation}
Scaling $z=z_++u/K$, the contour $\mathcal{C}$ may be approximated by the contour $\mathcal{E}$, which starts at $-\infty-i0$, below the real axis, circles the origin in the counterclockwise direction and ends at $-\infty+i0$, above the real axis. Thus, (\ref{S3_delta_rho>1}) becomes
\begin{equation}\label{S3_delta_rho>1_2}
\Delta H_K\sim K^{\alpha-1}z_+^{K-1}(z_+-z_-)^{\alpha-1}\Big[I_1+I_2/K+I_3/K^2+O(K^{-3})\Big],
\end{equation}
where the $I_j$ are expressible in terms of the $\Gamma$ function, with
\begin{eqnarray*}
I_1&\equiv&\frac{1}{2\pi i}\int_{\mathcal{E}}e^{u/z_+}(z_+-1)u^{-\alpha}du\\
&=&z_+^{1-\alpha}(z_+-1)/\Gamma(\alpha),
\end{eqnarray*}
and $I_2$ and $I_3$ can be given as similar contour integrals.
To find the dominant singularity, (\ref{S3_delta_rho>1_2}) must vanish. Viewing the $I_j=I_j(\theta)$ as functions of $\theta$, $I_1$ will vanish when $z_+=z_+(\theta)=1$, which occurs at $\theta=0$ if $\rho>1$. Thus we expand $\theta$ as $\theta=\alpha_2K^{-1}+\beta_2K^{-2}+O(K^{-3})$. Then (\ref{S3_delta_rho>1_2}) yields 
$$I_1'(0)\theta+\frac{1}{2}I_1''(0)\theta^2+K^{-1}\big[I_2(0)+I_2'(0)\theta\big]+K^{-2}I_3(0)=O(K^{-3}),$$
so that $\alpha_2=-I_2(0)/I_1'(0)$ and
$$\beta_2=-\Big[\frac{1}{2}I_1''(0)\alpha_2^2+I_2'(0)\alpha_2+I_3(0)\Big]\Big/I_1'(0),$$
 which leads to $\theta_s\sim-1/K-1/[(\rho-1)K^2]$. We note that to obtain the higher order approximations in (\ref{S2_th2_rho>1}), we needed to compute two more terms in (\ref{S3_delta_rho>1_2}). This concludes our derivation.

\end{document}